\documentclass[article, 12pt]{aiaa-pretty}    

\usepackage{amssymb}
\usepackage{amsmath}
\usepackage{booktabs}
 \usepackage{varioref}
 \usepackage{wrapfig}
 \usepackage{threeparttable}
 \usepackage{dcolumn}
  \newcolumntype{d}{D{.}{.}{-1}}
 \usepackage{nomencl}
  \makeglossary
 \usepackage{subfigure}
 \usepackage{graphicx}
 \usepackage{subfigmat}
 \usepackage{fancyvrb}
  \fvset{fontsize=\footnotesize,xleftmargin=2em}
 \usepackage{lettrine}
\usepackage{cite}
\usepackage{arydshln}
\usepackage{letltxmacro}
\LetLtxMacro{\originaleqref}{\eqref}
\renewcommand{\eqref}{Eq.~\originaleqref}

\newtheorem{theorem}{Theorem}

\newtheorem{remark}{Remark}
\newtheorem{lemma}{Lemma}
\newtheorem{definition}{Definition}

\newenvironment{proof}{{\it Proof. }}{\hfill $\Box$}

\newcommand{\Real}{\mathbb R}
\newcommand{\set}[1]{\left\{#1\right\}}
\newcommand{\real}[1]{{\mathbb R}^{#1}}
\newcommand{\bb}{{\boldsymbol b}}

\newcommand{\bd}{{\boldsymbol d}}

\newcommand{\bg}{{\boldsymbol g}}

\newcommand{\bw}{{\boldsymbol w}}

\newcommand{\by}{{\boldsymbol y}}

\newcommand{\bA}{{\boldsymbol A}}
\newcommand{\bB}{{\boldsymbol B}}

\newcommand{\bzero}{{\bf 0}}

\newcommand{\bpsi}{{\mbox{\boldmath $\psi$}}}

\newcommand{\bxi}{{\mbox{\boldmath $\xi$}}}

\newcommand{\bchi}{\mbox{\boldmath$\chi$}}

\author{ %
Isaac M. Ross\thanks{Distinguished Professor and Program Director, Control and Optimization, Department of Mechanical and Aerospace Engineering. imross@nps.edu }\\
\textit{Naval Postgraduate School, Monterey, CA 93943}
}
\title{Hessians in Birkhoff-Theoretic Trajectory Optimization}

\begin{document}
\maketitle


\section{Introduction}\label{sec:Intro}
Deep inside a trajectory optimization algorithm is an iterative scheme given by,
\begin{equation}\label{eq:iteration}
\bchi^{k+1} = \bchi^k + \delta\bchi^k, \quad k = 0, 1, \ldots
\end{equation}
where $\bchi^k \in \real{n}$ is an optimization variable and $\delta\bchi^k \in \real{n}$ is a search vector. In a line-search strategy\cite{NW:NumOptBook,luenberger-2008,bazaraa-2006}, the search vector is written as $\delta\bchi^k = \alpha^k \bd^k, \alpha^k \in \Real, \bd^k \in \real{n}$, where $\alpha^k > 0$ is called the step length and $\bd^k$ is called the search direction. In this approach, the search direction $\bd^k$ is found first and the step length $\alpha^k$ is computed afterwards. In a trust region strategy\cite{luenberger-2008,NW:NumOptBook,bazaraa-2006}, the search vector is found by constraining it to a ``trust'' region $\norm{\delta\bchi^k} \le \Delta^k > 0$  where $\norm{\cdot}$ is typically a $2$-norm and $\Delta^k$ is the trust-region radius. This approach can be contrasted to the line-search strategy as finding the ``step length'' $\Delta^k$ first and computing the search direction $\delta\bchi^k$ afterwards. In either case, the computation of $\delta\bchi^k$ involves solving a linear matrix equation,
\begin{equation}\label{eq:Ax=b}
\bA^k \bxi^k = \bg^k,  \quad k = 0, 1, \ldots
\end{equation}
where $\bA^k$ is an $n \times n$ matrix that contains an amalgamation of the discretized Jacobians and Hessians associated with the trajectory optimization problem, $\bg^k \in \real{n}$ comprises function and Jacobian evaluations based on the problem data functions and $\bxi^k$ is either $\bd^k$, as in the case of a line-search strategy or $\delta\bchi^k$ if the algorithm employs a trust-region.  The specifics of the process involved in generating \eqref{eq:Ax=b} from \eqref{eq:iteration} is described in detail in Section~\ref{sec:HamilFormOfHess}.

Solving \eqref{eq:Ax=b} for $\bxi^k$ is one of the most computationally expensive components of a trajectory optimization algorithm.  It can be roughly stated as having an asymptotic worst-case computational complexity\cite{trefethen-bau-1997,nesterov-book-2004,golub-book-2013} of $\mathcal{O}(n^3)$.  For example, Gauss elimination requires $(2/3)n^3$ floating point operations\cite{trefethen-bau-1997,golub-book-2013}. The $\mathcal{O}(n^3)$-statement is a more precise version of the well-known adage that if $n$ is small/large the computational time for solving \eqref{eq:Ax=b} is correspondingly small/large.  Furthermore, because \eqref{eq:Ax=b} must be solved for every $k = 0, 1, \ldots$, the worst-case computational time is additive and can be expressed as $\mathcal{O}(K n^3)$, where $K$ is the number of iterations.  If $\bA^k$ is structured, the computational complexity for solving \eqref{eq:Ax=b} can be reduced by employing non-generic methods and algorithms for solving linear matrix equations\cite{trefethen-bau-1997}. A trivial example of the last statement is the case when $\bA^k$ is diagonal. In this case, the computational complexity of solving \eqref{eq:Ax=b} is just $\mathcal{O}(n)$. To briefly put the arguments of the preceding statements in perspective, consider a value of $n=1,000$, a small number for a present-day computer. The computational time for $\mathcal{O}(n^3)$ operations is proportional to $10^9$ time units.  In contrast, an $\mathcal{O}(n)$ method takes a mere $10^3$ time units.  Stated differently, the $\mathcal{O}(n)$-computation is \emph{a million-times faster} than $\mathcal{O}(n^3)$ (for $n=1,000$).
In view of all these considerations, it is apparent that understanding the detailed structure of $\bA^k$ is an important problem for generating fast, stable and accurate solutions to trajectory optimization problems.

In this note, we derive the matrix $\bA^k$ for the universal Birkhoff-theoretic method for trajectory optimization developed in \cite{newBirk-part-I,newBirk-part-II,newBirk-2023}.  The contributions of this note are as follows:
\begin{enumerate}
\item We use an unconventional weighted Lagrangian and a special scaling of the optimization variables used in \cite{newBirk-part-I} to derive $\bA^k$. We show that $\bA^k$ can be split into two blocks: one that contains only the Birkhoff-specific components and another that depends only on the problem-specific data functions.  That is, the Birkhoff-specific components are decoupled from the problem-specific elements.
\item We show that the problem-specific components of $\bA^k$ are block diagonal matrices that can be expressed explicitly and only pointwise in terms of the data functions (and their derivatives).
\item We show that a particular arrangement of alternating primal and dual variables together with a particular listing of the constraints generates a matrix $\bA^k$ whose diagonal entries are mostly identity matrices.
\item Besides expressing it in an interesting manner, we show that the off-diagonal blocks of $\bA^k$ can be directly connected, via an identity transform, to the Hessian of the Hamiltonian, the Hessian of the endpoint Lagrangian and the Jacobian of the problem data functions.
\item To demonstrate the possible utilization of the special structure of $\bA^k$, we develop an eigenvalue theorem for the Birkhoff Hessian based on the well-known Gershgorin theorem\cite{kreyszig-2011}. Key results of this new theorem are:
    \begin{enumerate}
    \item The spectral radius of the Birkhoff Hessian is mesh-independent;
    \item Approximately,  $80\%$ of the eigenvalues are independent of the problem data; and,
    \item The real components of the problem-independent eigenvalues (i.e., approximately $80\%$ of the eigenvalues) are all contained in the narrow interval $[-2, 4]$.
    \end{enumerate}
\item We use the preceding results to quantify the computational complexity\cite{complexity-book-2007} of solving \eqref{eq:Ax=b} and point a way forward on how to address the grand challenge of solving a million-point trajectory optimization problem\cite{ross-million-2017}.
\end{enumerate}

The results of this note also suggest that it is possible to make informed changes to the trajectory optimization problem (such as coordinate transformations, scaling etc.\cite{scaling}) with immediate and detailed knowledge of their impact on the  innermost component of an algorithm, namely  \eqref{eq:iteration} and \eqref{eq:Ax=b}.

\section{An Overview of Birkhoff-Theoretic Discretizations of Optimal Control Problems}\label{sec:funda}

\begin{figure}[h!]
      \centering
      {\parbox{\columnwidth}{
      \centering
      {\includegraphics[width = 0.7\columnwidth]{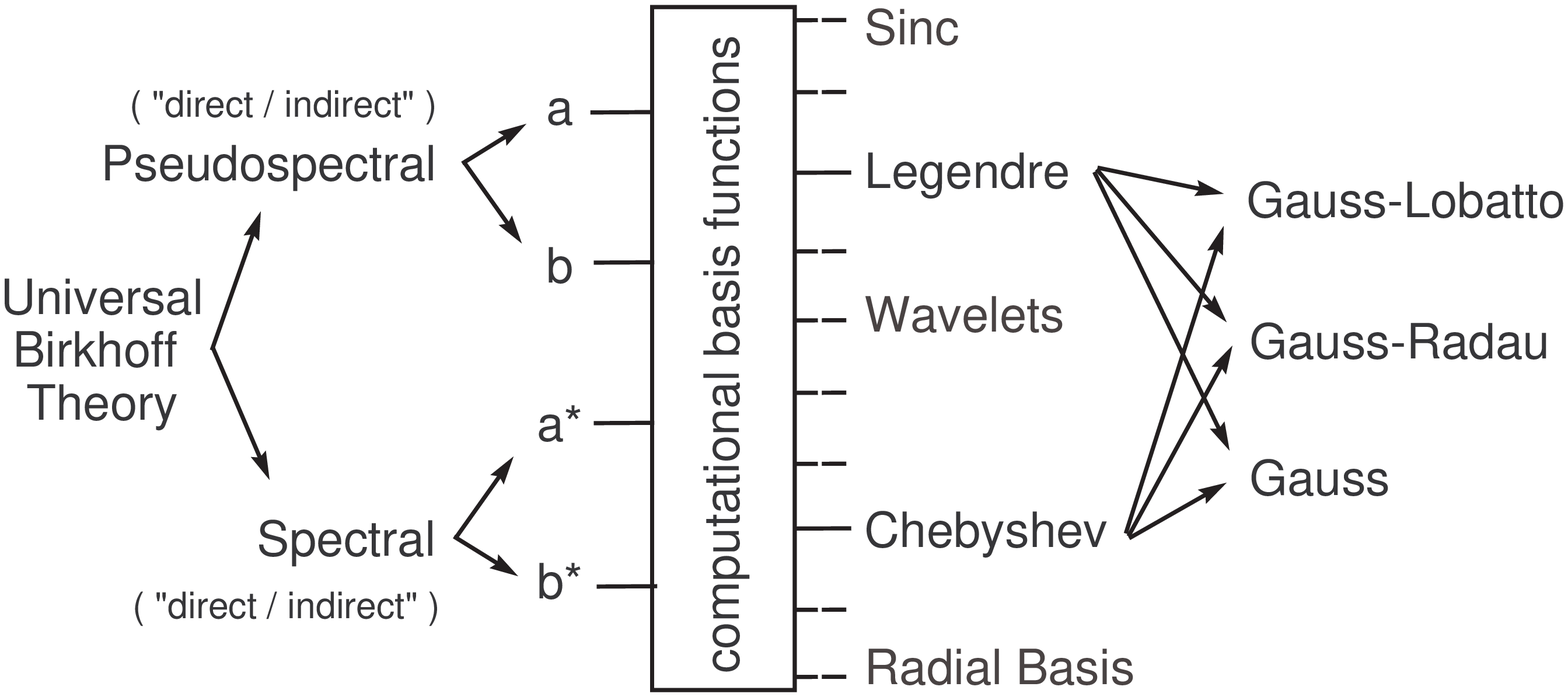}}
      \caption{\textsf{Schematic for illustrating the plethora of possible Birkhoff-theoretic methods for trajectory optimization; figure adapted from \cite{newBirk-part-II} and\cite{DIDO:arXiv}.}}\label{fig:BirkhoffPStypes}
      }
      }
\end{figure}

As illustrated in Fig.~\ref{fig:BirkhoffPStypes}, a universal Birkhoff-theoretic method for trajectory optimization\cite{newBirk-part-I,newBirk-part-II,DIDO:arXiv} generates an extremely large number of variants.
%
%
To limit the scope of this note, we consider only those Birkhoff methods that use orthogonal polynomials as basis functions.  In particular, we limit the scope of this note to the Birkhoff theory that utilizes only the Chebyshev and Legendre polynomial basis functions\cite{newBirk-part-I,newBirk-part-II,newBirk-2023}.  Because of the choice of these Gegenbauer (i.e., Chebyshev and Legendre) basis functions, we further limit the scope of this note by constraining the independent variable to the domain $[-1, 1] = [\tau^a, \tau^b]$ in order to circumvent the impact of time domain transformation techniques\cite{ross-book,advances,Radau-JGCD}.  Finally, as has been noted many times before\cite{fastmesh,furtherResults,advances,newBirk-part-I}, it is more informative to explain the key ideas by formulating a distilled optimal control problem that avoids the distractions of bookkeeping the multitude of variables and constraints.  Because the main point of this note is to illuminate the structure of the Hessian, we define a particular version of a distilled optimal control problem as follows:
%
\begin{eqnarray}
& (P) \left\{
\begin{array}{lrl}
\emph{Minimize } && J[x(\cdot), u(\cdot)] := E(x(\tau^a), x(\tau^b))\\[0.5em]
\emph{Subject to}&& \dot x(\tau) = f(x(\tau), u(\tau)) \\
                 && e(x(\tau^a), x(\tau^b)) =  0
\end{array} \right. & \label{eq:probP}
\end{eqnarray}
The data in Problem~$(P)$, as indicated in \eqref{eq:probP}, are as follows:
$J$ is a cost functional given by a twice-differentiable endpoint (``Mayer''\cite{ross-book}) cost function $E: \Real \times \Real \to \Real$; the pair, $\big(x(\cdot), u(\cdot)\big)$, is the unknown system trajectory (i.e., state-control function pair, $\tau \mapsto \big(x(t), u(t)\big) \in \Real \times \Real$); $f: \Real \times \Real \to \Real$ is a given twice-differentiable dynamics function; and, $e: \Real \times \Real \to \Real$ is a given twice-differentiable endpoint constraint function that constrains the endpoint values of $x(\tau)$ at $\tau^a = -1$ and $\tau^b = 1$.

A universal Birkhoff method\cite{newBirk-2023,newBirk-part-I} begins by choosing a grid $\pi^N := \set{\tau_0, \ldots, \tau_N},  \tau^a \le \tau_0 < \cdots < \tau_N \le \tau^b $ that allows $\tau \mapsto x(\tau)$ to be approximated by two equivalent $a$- and $b$-expansions:
\begin{subequations}\label{eq:xN=Birka+b}
\begin{align}
x^N(\tau) := x^a B_0^0(\tau) + \sum_{j=0}^N v_j B_j^a(\tau)\label{eq:xN=Birk-a-1D}\\
x^N(\tau) := \sum_{j=0}^N v_j B_j^b(\tau) + x^b B_N^N(\tau) \label{eq:xN=Birk-b-1D}
\end{align}
\end{subequations}
where $(x^a, x^b) = (x(\tau^a), x(\tau^b))$, $B_0^0(\tau)$, $B^N_N(\tau)$, $ B_j^a(\tau)$ and $B_j^b(\tau)$, $j = 0, \ldots N$ are the Birkhoff basis functions that satisfy the following interpolating conditions\cite{newBirk-2023,newBirk-part-I}:
\begin{align}\label{eq:birk-conditions}
\begin{aligned}
B_0^0(\tau^a)           &= 1,            &&    B_N^N(\tau^b)   = 1,                      \\
 {d_\tau B}_0^0(\tau_i)   &=  0,          &&   {d_\tau B}_N^N(\tau_i)   = 0,   & i = 0, \ldots, N \\
B_j^a(\tau^a)           &=  0,            && B_j^b(\tau^b) = 0,                       & j = 0, \ldots, N \\
{d_\tau B}_j^a(\tau_i) &=  \delta_{ij},  && {d_\tau B}_j^b(\tau_i) =  \delta_{ij}, &   i = 0, \ldots, N,  & j = 0, \ldots, N
\end{aligned}
\end{align}
In \eqref{eq:birk-conditions}, $d_\tau$ is the derivative operator, $d/d\tau$, and  $\delta_{ij}$ is the Kronecker delta.  The variables $v_j, j = 0, \ldots, N$ in \eqref{eq:xN=Birka+b} are called the virtual variables or the virtual control variables\cite{DIDO:arXiv,newBirk-part-I}.

Let, $x^N(\tau_i) = x_i$ and $u(\tau_i) = u_i$, $i = 0, \ldots, N$.  Define\cite{newBirk-2023,newBirk-part-I},
\begin{equation}\label{eq:manyDefs}
\begin{aligned}
X := & [x_0, \ldots, x_N]^T \\
V := & [v_0, \ldots, v_N]^T  \\
U := & [u_0, \ldots, u_N]^T \\
f(X, U) := & [f(x_0, u_0), \ldots, f(x_N, u_N)]^T  \\
\bw_B = & [w_0, \ldots, w_N]^T\\
\bb:= & [1, \ldots, 1]^T \\
[\bB^a]_{ij} := & [B^a_j({\tau_i})]  \\
[\bB^b]_{ij} := & [B^b_j({\tau_i})]  \\
N_n := & N+1
\end{aligned}
\end{equation}
where $w_i, i = 0, \ldots, N$ are the Birkhoff quadrature weights\cite{newBirk-2023}.  Then, the $a$-version of the Birkhoff discretization of Problem~$(P)$ is given by\cite{newBirk-part-I},
%
\begin{align}
X \in \real{N_n}, \quad U \in \real{N_n}, \quad V \in \real{N_n}, \quad (x^a, x^b) \in \real{2}   & \nonumber\\
 \left(\textsf{$P_a^N$}\right) \left\{
\begin{array}{lrl}
\emph{Minimize } & J^N_a[X, U, V, x^a, x^b] :=& E(x^a,x^b)\\
\emph{Subject to} & X  = &x^a\,\bb + \bB^a V \\
& V =& f(X, U)  \\
& x^b = & x^a + \bw_B^T V \\
& e(x^a, x^b)  = & 0
\end{array} \right. & \label{eq:Prob-PNa}
\end{align}
%
Problem~$(P^N_a)$ is a discretized optimal control problem (DOCP).  The discretization refers to time while the optimization variables $X, V, U, x^a$ and $x^b$ are all continuous.  In much of the trajectory optimization literature\cite{conway:survey,trelat:survey}, DOCPs (such as Problem~$(P^N_a)$) are called nonlinear programming problems (NLPs).  In principle, this is true; however, note that the variables $X, V, U, x^a$ and $x^b$ have different and special properties.  In an NLP these variables are treated as a single parameter, $\by := (X, V, U, x^a, x^b)$ that is largely agnostic to any special properties of its constituent variables vis-a-vis its continuous-time counterparts.  As noted elsewhere\cite{DIDO:arXiv,LRO-CSM,ross:shooting-arXiv},  DOCPs occupy a much smaller ``space'' than the ``space of all NLPs.'' In view of this, we approach a DOCP as a very special continuous optimization problem that is linked to Problem~$(P)$ in every step of the resulting analysis up to and including a connection to \eqref{eq:iteration} and \eqref{eq:Ax=b}. It will be apparent in the sections to follow that this unique treatment of Problem~$(P^N_a)$ is sharply different from the conventional ``NLP blackbox'' perspective.

\section{Connecting a Special Lagrangian to the Pontryagin Hamiltonian}\label{sec:Lag2Ham}
To draw a connection between a Lagrangian of a DOCP and the Pontryagin Hamiltonian of Problem~$(P)$, we need the following definition and lemma (proved in \cite{newBirk-part-I}):
\begin{definition}
The Birkhoff quadrature matrix is the following $N_n \times N_n$ matrix,
\begin{equation}
W_B:= diag[\bw_B] :=  \left[
\begin{array}{ccc}
w_0 & 0 & 0 \\
0 & \ddots & 0 \\
0 & 0 & w_N \\
\end{array}
\right]
\end{equation}
where $\bw_B$ is as given in \eqref{eq:manyDefs}.
\end{definition}
\begin{lemma}[\cite{newBirk-part-I}]\label{lemma:ross-1}
Given any $\epsilon > 0$, there exists an $N_\epsilon \in \mathbb{N}$ such that for all $N \ge N_\epsilon$ the following holds:
$$ W_B \bB^b + [\bB^a]^T W_B = \bzero  $$
\end{lemma}
Based on results of Theorem~3 of \cite{newBirk-part-I}, we scale the variables $X, V$ and $U$ in the DOCP given by Problem~$(P^N_a)$ by $W_B$ such that the new optimization variables $\widetilde{X}, \widetilde{V}$ and $\widetilde{U}$ are given by,
\begin{align*}
\widetilde{X} &= \bw_B \circ X = W_B X\\
\widetilde{V} &= \bw_B \circ V = W_B V \\
\widetilde{U} &= \bw_B \circ U = W_B U
\end{align*}
where $\circ$ denotes a Hadamard operation.  In not treating a DOCP as a generic NLP, consider the following specially weighted Lagrangian for Problem~$(P^N_a)$:
\begin{multline}\label{eq:Lagrangian}
L(\Omega, \Lambda, \lambda^b, \nu, \widetilde{X}, \widetilde{U}, \widetilde{V}, x^a, x^b)
 := E(x^a, x^b) \\
+ \Omega^T W_B \left(W_B^{-1} \widetilde{X} - \bB^a W_B^{-1} \widetilde{V} - x^a\bb \right)\\
+ \Lambda^T W_B \left(f(W_B^{-1} \widetilde{X},\, W_B^{-1} \widetilde{U}) - W_B^{-1} \widetilde{V}\right) \\
+ \lambda^b\left(x^a + \bw_B^T W_B^{-1} \widetilde{V} - x^b \right)
+ \nu e(x^a, x^b)
\end{multline}
where $\Omega$, $\Lambda$, $\lambda^b$ and $\nu$ are the Lagrange multipliers associated with the constraints implied in \eqref{eq:Lagrangian}.
Under these constructs, the necessary conditions for the variable-scaled DOCP are given by,
\begin{equation}\label{eq:gradL=0-pre}
\partial_{\widetilde{X}} L = 0 =  \partial_{\widetilde{V}} L = \partial_{\widetilde{U}} L = \partial_{x^a} L = \partial_{x^b} L
\end{equation}
To apply \eqref{eq:gradL=0-pre} in an efficient manner, it is useful to note that
\begin{align*}
\partial_X f(X, U)  & =
\left[
\begin{array}{ccc}
\partial_xf(x_0, u_0) & 0 & 0 \\
0 & \ddots & 0 \\
0 & 0 & \partial_xf(x_N, u_N) \\
\end{array}
\right] \\
& := diag[\partial_xf(x_0, u_0), \ldots, \partial_xf(x_N, u_N)]^T
\end{align*}
Defining\cite{newBirk-part-I},
\begin{equation}\label{eq:fx(X,U)}
\partial_x f(X, U) := [\partial_xf(x_0, u_0), \ldots, \partial_xf(x_N, u_N)]^T
\end{equation}
it follows that we can write,
$$\partial_X f(X, U) = diag[\partial_x f(X,U)]$$
Hence we have,
$$ [\partial_X f(X, U)]^T\Lambda = diag[\partial_x f(X,U)] \Lambda =  \partial_x f(X,U) \circ \Lambda $$
With these notational conveniences, an application of \eqref{eq:gradL=0-pre} using \eqref{eq:Lagrangian} and Lemma~\ref{lemma:ross-1} generates the following set of equations,
\begin{subequations}\label{eq:gradL=0}
\begin{align}
 \Omega +  \partial_x f(X, U)\circ \Lambda = & 0 \label{eq:adjoint-Omega}\\
 \Lambda - \lambda^b\, \bb - \bB^b \Omega = & 0 \label{eq:adjoint-linear} \\
 \partial_u f(X, U)\circ \Lambda  = & 0 \label{eq:HMC-0}\\
\lambda^a + \partial_{x^a} \overline{E}(\nu, x^a, x^b) = & 0 \label{eq:tvc-a}\\
 \lambda^b - \partial_{x^b} \overline{E}(\nu, x^a, x^b) = & \label{eq:tvc-b}0
\end{align}
\end{subequations}
where $\partial_u f(X,U)$ is defined similarly to \eqref{eq:fx(X,U)},  $\lambda^a$ is defined by\cite{newBirk-part-I},
\begin{equation}\label{eq:lam-a-def}
\lambda^a := \lambda^b -  \bw_B^T \Omega
\end{equation}
and $\overline{E}(\nu, x^a, x^b)$  is the endpoint Lagrangian\cite{ross-book} function (associated with Problem~$(P)$) defined by,
\begin{equation}
\overline{E}(\nu, x^a, x^b) :=  E(x^a, x^b) + \nu\, e(x^a, x^b)
\end{equation}

The purpose of the preceding mathematical machinations is now clear, at least, partially. Observe that \eqref{eq:tvc-a} and \eqref{eq:tvc-b} are, in fact, the transversality conditions\cite{ross-book,brysonHo} for Problem~$(P_a)$.  In fact, the remainder of the equations in \eqref{eq:gradL=0} are discretizations of the remainder of the first-order necessary conditions of Pontryagin\cite{ross-book}.  To make the preceding statement more precise, let $H(\lambda, x, u)$ denote the Pontryagin Hamiltonian\cite{ross-book} for Problem~$(P)$ given by,
\begin{equation}
H(\lambda, x, u) := \lambda\, f(x, u)
\end{equation}
The necessary condition for Hamiltonian minimization is given by,
\begin{equation}\label{eq:Hu=0}
 \partial_u H(\lambda(\tau), x(\tau), u(\tau)) = \lambda(\tau)\, \partial_u f(x(\tau), u(\tau))  = 0  \quad \forall\ \tau
 \end{equation}
Hence, over the grid $\pi^N$, \eqref{eq:Hu=0} maps to,
\begin{equation}\label{eq:Hu=0-i}
 \partial_u H(\lambda(\tau_i), x(\tau_i), u(\tau_i)) = \lambda(\tau_i)\, \partial_u f(x(\tau_i), u(\tau_i ))  = 0  \qquad
 \forall\ \tau_i, \ i = 0, \ldots, N
 \end{equation}
Using the notational convention defined in  \eqref{eq:fx(X,U)}, it follows that \eqref{eq:Hu=0-i} can be written as,
\begin{equation}\label{eq:Hu-connect}
\partial_u H(\Lambda, X, U) = \Lambda\,\circ \partial_u f(X, U)  = 0
\end{equation}
Hence, \eqref{eq:HMC-0} is identical to the Hamiltonian form given by \eqref{eq:Hu-connect}.  Similarly, it is straightforward to show that
\begin{equation}\label{eq:Hx}
\partial_x H(\Lambda, X, U) = \Lambda\,\circ \partial_x f(X, U)
\end{equation}
Using \eqref{eq:Hx} it can be shown\cite{newBirk-part-I} that \eqref{eq:adjoint-Omega} and \eqref{eq:adjoint-linear} are the Birkhoff discretizations of the adjoint equation,
\begin{equation}
-\dot\lambda(\tau) =   \partial_x H(\lambda, x, u)
\end{equation}
written in terms of the differential-algebraic equation,
\begin{equation}
\dot\lambda(\tau) :=  \omega(\tau), \quad \omega(\tau) + \partial_x H(\lambda, x, u) = 0
\end{equation}
where $\omega(\tau)$ is called the co-virtual variable\cite{newBirk-part-I,DIDO:arXiv}.
Collecting all relevant equations, it follows that \eqref{eq:gradL=0} and hence \eqref{eq:gradL=0-pre} can be directly connected to Problem~$(P)$ in the sense that it can be rewritten in terms of the discretized Pontryagin Hamiltonian as follows:
\begin{subequations}\label{eq:gradL2H=0}
\begin{align}
 \Omega +  \partial_x H(\Lambda, X, U) = & 0 \\
 \Lambda - \lambda^b\, \bb - \bB^b \Omega = & 0  \\
 \lambda^a - \lambda^b +  \bw_B^T \Omega = & 0\\
 \partial_u H(\Lambda, X, U)  = & 0 \label{eq:Hu=0-strong}\\
\lambda^a + \partial_{x^a} \overline{E}(\nu, x^a, x^b) = & 0\\
 \lambda^b - \partial_{x^b} \overline{E}(\nu, x^a, x^b) = & 0
\end{align}
\end{subequations}
%

\begin{remark}\label{rem:LagLam}
When a generic unweighted Lagrangian is differentiated with respect to its multipliers, it regenerates the equality constraints of an NLP\cite{luenberger-2008,NW:NumOptBook,bazaraa-2006}.  It is straightforward to show that differentiating the special Lagrangian of \eqref{eq:Lagrangian} with respect to $\Omega$ and $\Lambda$ generates the constraints in Problem~$(P^N_a)$ in its spectral form\cite{newBirk-part-I}.  That is, the resulting constraints are multiplied by $W_B$. This is a weak representation of the constraints because for $N>2$, $0 < w_i < 1, i = 0, \ldots N$.  To impose the primal constraints strongly, we use its formulation in its original ``strong'' (pseudospectral) form given in \eqref{eq:Prob-PNa}.
\end{remark}
\begin{remark}\label{rem:lam-a}
Equation~(\ref{eq:lam-a-def}) is an ``artifice'' that does not naturally occur if Problem~$(P^N_a)$ is treated as a generic NLP even if the special Lagrangian of \eqref{eq:Lagrangian} is used. In the context of an NLP, \eqref{eq:lam-a-def} is simply a clever grouping of the NLP multipliers to produce a ``new'' linear equation. In the context of a trajectory optimization problem, \eqref{eq:lam-a-def} is a critical equation in the sense that $\lambda^a$ is the initial costate.
\end{remark}

Remarks~\ref{rem:LagLam} and \ref{rem:lam-a} illustrate why it may be unwise to treat Problem~$(P^N_a)$ a generic NLP.  In fact, the following sections amplify this point.

\section{Hamiltonian Form of the Birkhoff Hessian}\label{sec:HamilFormOfHess}
In connecting the ideas introduced in Section~\ref{sec:Intro} to Section~\ref{sec:Lag2Ham}, let $\bchi$ tentatively denote the arguments of the special Lagrangian given by \eqref{eq:Lagrangian}.  Then, from the results of Section~\ref{sec:Lag2Ham}, it follows that a solution to Problem~$(P^N_a)$ can be obtained by solving a system of linear and nonlinear equations that can be expressed as $\boldmath{\mathcal{F}}(\bchi) = \bzero$, the details of which will be presented shortly. Newton's method for solving a system of nonlinear equations is given by\cite{luenberger-2008,NW:NumOptBook,bazaraa-2006},
\begin{equation}\label{eq:Newton-method}
\partial_{\bchi}\boldmath{\mathcal{F}}(\bchi^k) \bd^k = -\boldmath{\mathcal{F}}(\bchi^k)
\end{equation}
where $\bd^k:= \bchi^{k+1} - \bchi^k$ is the Newton step.  See also \eqref{eq:iteration} and \eqref{eq:Ax=b} for context.  Because the source of $\boldmath{\mathcal{F}}(\bchi)$ is the gradient of the Lagrangian (see \eqref{eq:gradL=0}), the Jacobian in \eqref{eq:Newton-method} is the Hessian of $L(\Omega, \Lambda, \lambda^b, \nu, \widetilde{X}, \widetilde{U}, \widetilde{V}, x^a, x^b)$, under the caveat of Remark~\ref{rem:LagLam}.  Hence, we refer to such a Jacobian as a Birkhoff Hessian. The task at hand is to develop a Birkhoff Hessian that can be directly connected to Problem~$(P)$.

To construct a Birkhoff Hessian that is revelatory with respect to its properties, we first organize all the variables in a particular arrangement of mostly alternating primal and dual variables according to (see also Remark~\ref{rem:lam-a}),
\begin{multline}\label{eq:chi:=}
\bchi:= (X, \Lambda, V,  \Omega, U, x^a, \lambda^b, x^b, \nu, \lambda^a)  \\
\in \real{N_n}\times \real{N_n} \times \real{N_n} \times \real{N_n} \times \real{N_n} \times \Real \times \Real \times \Real \times \Real \times \Real
\end{multline}
For similar reasons, we also arrange the components of $\boldmath{\mathcal{F}}(\bchi)$ (i.e., \eqref{eq:gradL2H=0} and the constraints in \eqref{eq:Prob-PNa}) in the following specific manner:
\begin{subequations}
\begin{align}
\mathcal{F}_1(\bchi) &= X - x^a\,\bb - \bB^a V         \\
\mathcal{F}_2(\bchi) &= \Lambda - \lambda^b\, \bb - \bB^b \Omega      \\
\mathcal{F}_3(\bchi) &= V - \partial_\lambda H(\Lambda, X, U)     \\
\mathcal{F}_4(\bchi) &= \Omega +  \partial_x H(\Lambda, X, U)    \\
\mathcal{F}_5(\bchi)&= \partial_u H(\Lambda, X, U)  \\
\mathcal{F}_6(\bchi) &=  x^a -  x^b + \bw_B^T V       \\
\mathcal{F}_7(\bchi)&= -\lambda^a +  \lambda^b - \bw_B^T \Omega       \\
\mathcal{F}_8(\bchi) &= -\lambda^b + \partial_{x^b} \overline{E}(\nu, x^a, x^b)     \\
\mathcal{F}_{9}(\bchi)&= \partial_{\nu} \overline{E}(\nu, x^a, x^b)  \\
\mathcal{F}_{10}(\bchi) &= \lambda^a + \partial_{x^a} \overline{E}(\nu, x^a, x^b)
\end{align}
\end{subequations}
Note the sign changes in $\mathcal{F}_7(\bchi)$ and $\mathcal{F}_8(\bchi)$.  Note also that $\mathcal{F}_3(\bchi)$ and $\mathcal{F}_{9}(\bchi)$ are the primal constraint equations (i.e., the constraints in \eqref{eq:Prob-PNa}) rewritten in terms of the Hamiltonian and the endpoint Lagrangian respectively using the fact that,
\begin{equation}
\partial_\lambda H(\lambda, x, u) = f(x,u), \quad \partial_{\nu} \overline{E}(\nu, x^a, x^b) = e(x^a, x^b)
\end{equation}
Hence (see \eqref{eq:fx(X,U)}),
\begin{align*}
\partial_\lambda H(\Lambda, X, U) & = f(X,U)  \\
 & :=  [f(x_0, u_0), \ldots, f(x_N, u_N)]^T \\
 & = \partial_\Lambda H(\Lambda, X, U)
\end{align*}
Furthermore,
\begin{equation}\label{eq:HessHam=diag}
\partial^2_{\Lambda, X} H(\Lambda, X, U) = \partial_X f(X, U) 
=
\left[
\begin{array}{ccc}
\partial_xf(x_0, u_0) & 0 & 0 \\
0 & \ddots & 0 \\
0 & 0 & \partial_xf(x_N, u_N) \\
\end{array}
\right]
\end{equation}
Similarly, it is straightforward to show that the second derivative of the Hamiltonian with respect to any of its variables $(\Lambda, X, U)$ is a diagonal matrix.

To simplify the identification of the matrix subblocks  in the construction of the Birkhoff Hessian, we use the following notational convenience:
$$ \partial^2_{\Lambda, X} H(\Lambda, X, U) \equiv H_{\Lambda, X}, \quad  \partial^2_{x^a, x^b}\overline{E} \equiv \overline{E}_{a, b} $$
With these conveniences, the Birkhoff Hessian associated with the DOCP given by Problem~$(P^N_a)$ (i.e., Jacobian of $\boldmath{\mathcal{F}}(\bchi)$) can now be written succinctly as,
\begin{equation}\label{eq:Hessian=}
\bA := \left[
             \begin{array}{ccccc;{1pt/2pt}ccccc}
               I & 0 & - \bB^a & 0  & 0 & -\bb & 0 & 0 & 0 & 0\\
               0 & I & 0 & -\bB^b & 0 & 0 & -\bb & 0 & 0 & 0\\ 
               -H_{\Lambda, X}  & 0 & I  &  0 & -H_{\Lambda, U}  & 0 & 0 & 0 & 0 & 0\\
                H_{X,X}& H_{X,\Lambda} & 0  & I & H_{X, U} & 0 & 0 & 0 & 0 & 0\\
                H_{U,X} & H_{U,\Lambda} & 0 & 0 & H_{U, U} & 0 & 0 & 0 & 0 & 0 \\ \hdashline[1pt/2pt]
                0 & 0 & \bw_B^T & 0 & 0  & 1 & 0 & -1 & 0 & 0 \\
                0 & 0 & 0 & -\bw_B^T & 0 & 0 & 1 & 0 & 0 & -1\\
                0 & 0 & 0 & 0 & 0 &  \overline{E}_{b,a} &-1 & \overline{E}_{b,b} & \overline{E}_{b,\nu}& 0\\
                0 & 0 & 0 & 0 & 0 & \overline{E}_{\nu,a} &0 & \overline{E}_{\nu,b}  & 0 & 0\\
                0 & 0 & 0 & 0 & 0 & \overline{E}_{a,a} & 0 & \overline{E}_{a,b} & \overline{E}_{a,\nu} & 1 \\
             \end{array}
           \right]
\end{equation}
Some comments about $\bA$ are in order:
\begin{enumerate}
\item $I$ is an $N_n \times N_n$ identity matrix.  The $0$'s in \eqref{eq:Hessian=} are matrices of order $N_n \times N_n$ or $N_n \times 1$ or or $1 \times N_n$ or $1 \times 1$.  The precise size of these $0$'s are apparent from the context.
\item The Hessian is partitioned as indicated in \eqref{eq:Hessian=} to indicate that the top left ``$11$-block'' corresponds to continuous variables that have values over the entire grid while the bottom right ``$22$-block'' corresponds to continuous variables that take values only over the boundary points $a$ and $b$.
\item When compared to the $22$-block, the  $11$-block is relatively large.  The $22$-block is $5 \times 5$ while the $11$-block is $5N_n \times 5N_n$.  In other words, the $11$-block overwhelmingly dominates in size.
\item The variables and functions (i.e., $\bchi$ and $\boldmath{\mathcal{F}}$)  are arranged to show that the identity matrix occupies all the block diagonal elements of the Hessian in the top left block except for the lower right corner.  Similarly the diagonal elements of the bottom right block are $0$'s and $1's$ except for the the middle $33$-element of this sub-block.
\item The first two ``rows'' (block-wise) of the Hessian matrix, along with ``rows'' $6$ and $7$ are problem-independent. All other rows depend only the problem data functions and are decoupled from any specificity with respect to a Birkhoff-theoretic method for trajectory optimization.  Except for the bottom right block of the Hessian, all problem-dependent data are located only along the diagonals of the block elements indicated in \eqref{eq:Hessian=}.
\end{enumerate}

If Problem~$(P^N_a)$ were treated as a generic NLP, its Hessian would not be the same as \eqref{eq:Hessian=}.  This is because of the following reasons:
\begin{enumerate}
\item In a typical NLP, the special Lagrangian given by \eqref{eq:Lagrangian} is not used. The use of such a Lagrangian was critical to the development of \eqref{eq:Hessian=}.
\item Even if an NLP were to be modified to incorporate the special Lagrangian of \eqref{eq:Lagrangian}, \eqref{eq:lam-a-def} would not appear as part of its multiplier set. In this case the Hessian would be a ``$9 \times 9$'' matrix instead of the ``$10 \times 10$'' form presented in \eqref{eq:Hessian=}.  See Section~\ref{sec:NLPHessians} for additional details.
\item Because an NLP organizes variables differently than \eqref{eq:chi:=}, the identity matrices indicated in \eqref{eq:Hessian=} would lose their prime location along the diagonals.  See Section~\ref{sec:NLPHessians} for additional details.
\item The Hessian in an NLP is symmetric\cite{NW:NumOptBook,luenberger-2008,bazaraa-2006}.  It is obvious that the Birkhoff Hessian as expressed in \eqref{eq:Hessian=} is not symmetric.
\end{enumerate}
\begin{remark}
It is straightforward but tedious to show that no rearrangement of the rows/columns of $\bA$ will generate a symmetric matrix. A clue for this ``theorem of impossibility'' is the last column of $\bA$. The origin of this column (and $\mathcal{F}_7(\bchi)$) is the treatment of \eqref{eq:lam-a-def} as a constraint.
\end{remark}
The location of the identity matrices along the diagonal can be used judiciously in several different ways.  Complete details of utilizing this property for solving an optimal control problem are well beyond the scope of this note.  Nonetheless, we prove in the next section a general theorem pertaining to the eigenvalues of $\bA$.

\section{A Birkhoff Hessian Eigenvalue Theorem}
The eigenvalues of a matrix are a critically important component of its fundamental properties\cite{kreyszig-2011}. For a non-numerical matrix such as the one given by \eqref{eq:Hessian=}, it is typically fruitless to find its eigenvalues in terms of its constituent elements. However, Gershgorin's theorem\cite{kreyszig-2011} provides a simple and elegant summation formula to estimate the range of possible eigenvalues of a square matrix.  In this section, we prove an eigenvalue theorem by combining (a) Gershgorin's theorem, (b) the properties of the Birkhoff matrices $\bB^a$ and $\bB^b$\cite{newBirk-2023,newBirk-part-I} and (c) the structure of the Hessian as presented in \eqref{eq:Hessian=}.
%
\begin{theorem}\label{thm:ross-hess}
Let
\begin{enumerate}
\item $\pi^N$ be the Legendre/Chebyhsev (``Gegenbauer'') family of grids indicated in Fig.~\ref{fig:BirkhoffPStypes};
\item Problem~$(P)$, given by \eqref{eq:probP}, be discretized over $\pi^N$ by a universal Birkhoff-theoretic method; and
\item $\bA$, given by \eqref{eq:Hessian=}, be the Hessian matrix associated with a Birkhoff-Gegenbauer discretization of Problem~$(P)$.
\end{enumerate}
Then, the spectral radius of the Birkhoff Hessian is mesh/grid-independent; i.e., independent of $N$.  Furthermore, the real part(s) of:
\begin{enumerate}
\item $4N_n + 2$ eigenvalues are numerically bounded by the interval $[-2, 4]$, $(N_n = N+1)$.
\item One eigenvalue is bounded by the interval $[0, 2]$.
\item One eigenvalue is bounded by $\mp\left(\abs{\overline{E}_{\nu,a}} + \abs{\overline{E}_{\nu,b}}\right) $
\item One eigenvalue is bounded by $\abs{\overline{E}_{b,b}} \mp R^E$, where $R^E := 1 + \abs{\overline{E}_{b,a}} + \abs{\overline{E}_{b,\nu}} $.
\item $N_n$ eigenvalues are bounded by $\abs{H_{u_i, u_i}} \mp R^H_i, i = 0, \ldots, N$, where $R^H_i := \abs{H_{u_i, x_i}} + \abs{H_{u_i, \lambda_i}} $
\end{enumerate}
\end{theorem}
\begin{proof}
We will first prove Statements 1--5.  Then, the assertion that the spectral radius of the Hessian is independent of $N$ will be obvious.

We prove Statement 1 in three parts:
\begin{enumerate}
\item First, consider the identity matrices occupying the $11$- and $22$-blocks of $\bA$. For these blocks, consider applying Gershgorin's theorem row wise.  Each element of the vector $\bb$ is unity. Hence, the absolute value of any row of $-\bb$ is  unity. The sum of the absolute values of any row of $\bB^a$ or $\bB^b$ is less than or equal to $2$\cite{newBirk-2023}.  Hence the sum of the absolute values of the off diagonal elements of the first two rows of $\bA$ is $3$.  Applying Gershgorin's theorem, it follows that the real parts of $2N_n$ eigenvalues are centered at one with a radius of $3$.  Or equivalently, the real parts of $2N_n$ eigenvalues lie in the interval $[1-3, 1+3] = [-2, 4]$.
\item Next, consider the identity matrices occupying the $33$- and $44$-blocks of $\bA$.  For these blocks, we apply Gershgorin's theorem column wise. The absolute value of each element of $\bw_B$ is at most unity\cite{newBirk-2023}. Hence, it follows (from the results of the previous paragraph) that the real parts of an additional $2N_n$ eigenvalues lie in the interval $[-2, 4]$ raising the tally to $4N_n$.
\item Finally, consider the $1$'s occupying the $66$- and $77$ ``blocks'' of $\bA$.  Applying Gershgorin's theorem row wise to these blocks results in (see previous paragraph) two more eigenvalues lying in the interval $[-2, 4]$.
\end{enumerate}
Thus, Statement 1 is proved. To prove Statement 2 of the theorem, consider the last column of $\bA$. It is obvious that one eigenvalue is between $0$ and $2$.

To prove Statements 3, 4 and 5, we apply Gershgorin's theorem (row or column wise) to the diagonal entries of $\bA$ occupying the $99$-, $88$- and $55$-blocks of $\bA$ respectively.

From the proofs of Statements 1--5, it follows that the spectral radius of the Hessian is independent of $N$.
\end{proof}
\begin{remark}
For the $77$-block of $\bA$, applying Gershgorin's theorem along the column generate an estimate of an eigenvalue that is dependent upon $N$. Choosing row $7$ in the proof of Theorem~\ref{thm:ross-hess} was critical to proving mesh independence.
\end{remark}
\begin{remark}
Imbedded in  \eqref{eq:Hessian=} are at least six different possible Hessians based on the choice of the Gegenbauer grid (see Fig.~\ref{fig:BirkhoffPStypes}).  Hence Theorem~\ref{thm:ross-hess} imbeds at least six grid-dependent theorems.
\end{remark}

\section{A Preliminary Analysis of Computational Complexity}

The computational cost of solving solving \eqref{eq:Ax=b} with $\bA$ given by \eqref{eq:Hessian=} can be classified in terms of space and time complexity\cite{complexity-book-2007}.  To facilitate a simplified and practical discussion on this topic, a collection of space and time complexity numbers are provided in Table~\ref{table:size+speed} in a manner that is relatable to a present-day processor.
%
\begin{table}[h!]
\centering
\begin{tabular}{lcccccr}
  \toprule
 variables \qquad\qquad &\multicolumn{2}{c}{Memory (GB)} &\multicolumn{4}{c}{Compute Time (s/TFLOPS)}\\
  \cmidrule(lr){2-3} \cmidrule(lr){4-7}
  $n$ &  $n$ &  $n^2$  &  $n$  &$n\log(n)$ & $n^2$  & $n^3$\\
  \midrule
  1,000 & 0.000008 &0.008 & $1\times 10^{-9}$  & $7\times 10^{-9}$  & $0.000001$ & 0.001\\
  10,000 & 0.00008 &0.8  &  $1\times 10^{-8}$   & $9\times 10^{-8}$  & $0.0001$  & 1\\
  100,000 & 0.0008 & 80 & $1\times 10^{-7}$   &  $11 \times 10^{-7}$ & $0.01$ &  1,000\\
  1,000,000 & 0.008 & 8,000 & $1\times 10^{-6}$   &  $14\times 10^{-6}$ & $1$ & 1,000,000\\
  \bottomrule
\end{tabular}
\caption{Space and time complexity table baselined to gigabytes and terraflops}
\label{table:size+speed}
\end{table}
%
The computational space is computer memory listed in gigabytes (GB) in Table~\ref{table:size+speed} where $n$ is a generic number of variables.  Its relationship to the number of variables in trajectory optimization will be connected shortly. The computational time in Table~\ref{table:size+speed} is in seconds, baselined to a terraflops (TFLOPS) processor.  An example of how to read this table is as follows: Suppose we have a $1,000$-variable problem. Storing a square matrix for $n=1,000$ requires $0.008$ GB of memory ($=n^2$ in Table~\ref{table:size+speed}, namely, row 1, column 3).  If \eqref{eq:Ax=b} is solved using an algorithm that has $\mathcal{O}(n^3)$ complexity, then it can be solved in about $0.001$ seconds (row 1, column 7). We caveat this number and all of the discussions to follow with the fact that actual compute time depends on many more factors than FLOPS such as the number of pipelines, pipeline length, cache size, cache latency, instruction sets, number of registers, etc.\cite{computer-book-2019}.  Consequently, we use the compute times provided in Table~\ref{table:size+speed} to simply provide a relatable number in seconds. See also \cite{newBirk-part-I} for a critique on using actual computer time for solving a problem to indicate  ``fast'' in fast trajectory optimization.

As noted in Section~\ref{sec:Intro}, solving \eqref{eq:Ax=b} is one of the most computationally expensive components of a trajectory optimization algorithm. To analyze this computational cost, we permute the rows of $\bA$ given by \eqref{eq:Hessian=} to cast it in the following form:
\begin{equation}\label{eq:PA=}
P\bA := \left[
\begin{array}{c}
\bA_0\\ \hdashline[1pt/2pt]
\bA_{data}
\end{array}
\right]
:=
 \left[
             \begin{array}{cccccccccc}
               I & 0 & - \bB^a & 0  & 0 & -\bb & 0 & 0 & 0 & 0\\
               0 & I & 0 & -\bB^b & 0 & 0 & -\bb & 0 & 0 & 0\\ 
               0 & 0 & \bw_B^T & 0 & 0  & 1 & 0 & -1 & 0 & 0 \\
                0 & 0 & 0 & -\bw_B^T & 0 & 0 & 1 & 0 & 0 & -1\\ \hdashline[1pt/2pt]
               -H_{\Lambda, X}  & 0 & I  &  0 & -H_{\Lambda, U}  & 0 & 0 & 0 & 0 & 0\\
                H_{X,X}& H_{X,\Lambda} & 0  & I & H_{X, U} & 0 & 0 & 0 & 0 & 0\\
                H_{U,X} & H_{U,\Lambda} & 0 & 0 & H_{U, U} & 0 & 0 & 0 & 0 & 0 \\
                0 & 0 & 0 & 0 & 0 &  \overline{E}_{b,a} &-1 & \overline{E}_{b,b} & \overline{E}_{b,\nu}& 0\\
                0 & 0 & 0 & 0 & 0 & \overline{E}_{\nu,a} &0 & \overline{E}_{\nu,b}  & 0 & 0\\
                0 & 0 & 0 & 0 & 0 & \overline{E}_{a,a} & 0 & \overline{E}_{a,b} & \overline{E}_{a,\nu} & 1 \\
             \end{array}
           \right]
\end{equation}
In \eqref{eq:PA=}, $P$ is a permutation matrix that rearranges $\bA$ such that the ``top'' matrix $\bA_0$ is independent of the problem data and $\bA_{data}$ is the data-dependent ``bottom'' matrix.

\subsection{Computational Cost Associated with the Data-Dependent Matrix, $\bA_{data}$ }
\begin{figure}[h!]
      \centering
      {\parbox{\columnwidth}{
      \centering
      {\includegraphics[width = 0.65\columnwidth, trim=0in 0.4in 0in 0in, clip]{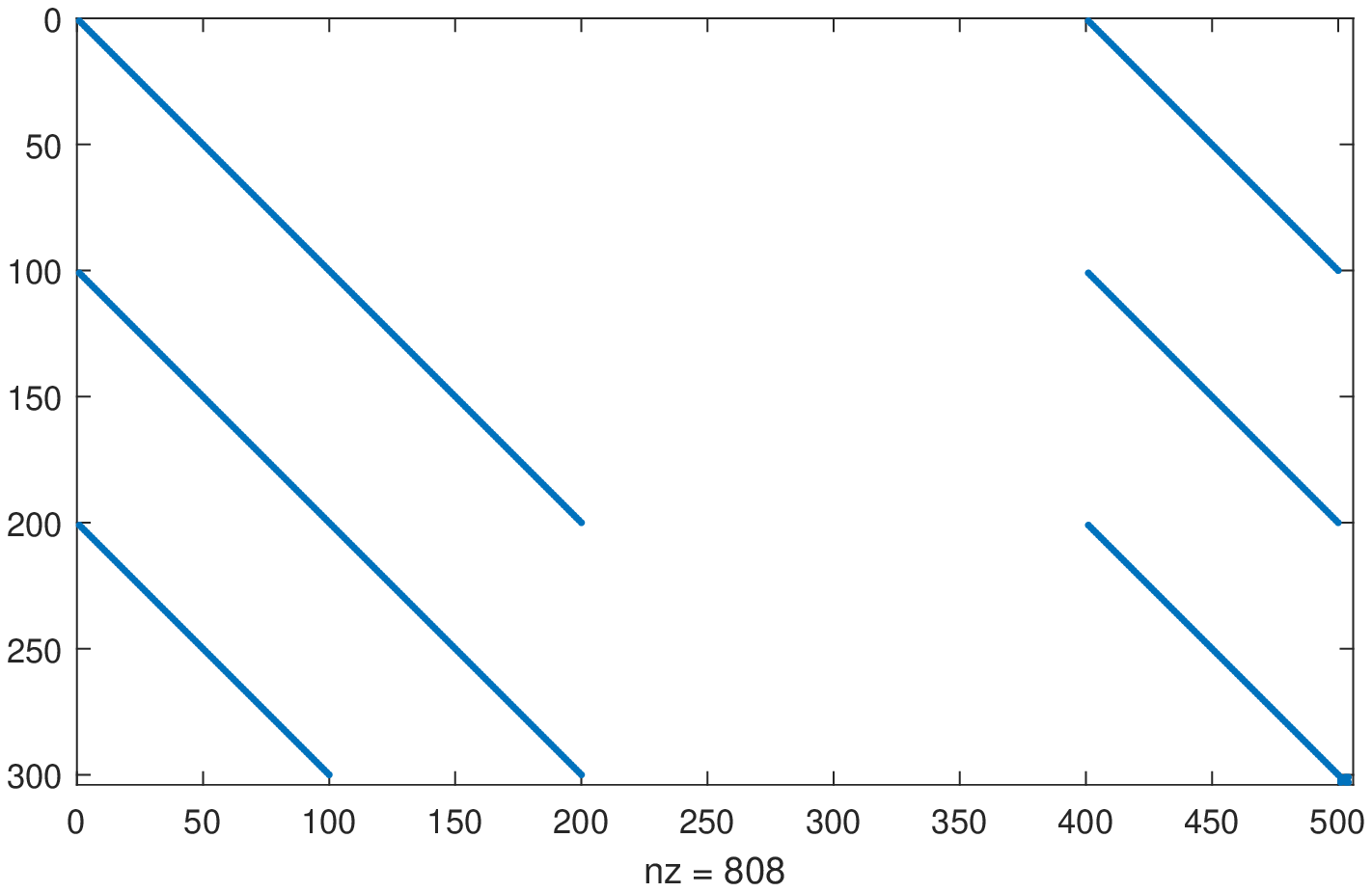}}
      \caption{\textsf{Data-dependent sparsity pattern of $\bA_{data}$}}\label{fig:Adata}
      }
      }
\end{figure}

From \eqref{eq:HessHam=diag} it follows that the Hessian of the Hamiltonian is a diagonal matrix (in its Birkhoff-discretized form).  Ignoring the constants (i.e., both ones and zeros) the sparsity pattern of $\bA_{data}$ is shown in Fig.~\ref{fig:Adata}.
%
%
It is apparent by inspection of this figure that space complexity (i.e. computer memory) is not an issue.  In terms of specifics, taking the symmetry of the Hessian of the Hamiltonian into account, the required memory for storing the Hessian is only $5N_n$ for Problem~$(P^N_a)$. In the general case, it is straightforward to show that the memory needs for $\bA_{data}$ is just $(2N_x^2 + 2N_x N_u +N_u^2) \times N_n $, where $N_x$ and $N_u$ are the number of state and control variables respectively (in continuous time). That is, the growth is only linear with respect to $N_n$. Hence, the Hessian storage requirements belongs to the second column of Table~\ref{table:size+speed}, namely $\mathcal{O}(n)$.  As an example, the computer memory for storing $\bA_{data}$ for a million grid points for a typical aerospace guidance problem (i.e. $N_x = 6, N_u = 3$) is less than one GB of memory ($\sim 0.9$ GB).

The computational cost of computing the product of $\bA_{data}$ with a vector (see \eqref{eq:Ax=b}) is $3N_n$ by inspection of Fig.~\ref{fig:Adata}.  In the general case, the computational cost is approximately $3N_x^2 \times N_n$.  Hence, the computational speed of computing a matrix-vector product with $\bA_{data}$ as the matrix is just $\mathcal{O}(N_n)$.  From the fourth column of  Table~\ref{table:size+speed}, it follows that this computation can be done in tens of microseconds even when the grid is one million points.

It is clear that neither space nor time complexity associated with $\bA_{data}$ is a present-day obstacle to the challenge posed in \cite{ross-million-2017} of producing a real-time, million-point trajectory optimization solver.  More importantly, because the contents of $\bA_{data}$ are problem specific, the preceding analysis is valid across the scope of all trajectory optimization problems.  The origin of this powerful conclusion is due to the following two facts: (i) A universal Birkhoff-theoretic method decouples the discretization from the problem-specific data and (ii) the problem-specific data is sampled pointwise with no dependencies on the adjacent grid points.  The last point is in sharp contrast to Runge-Kutta discretization methods that have interdependencies across the grid points\cite{ascher-petzold-1998}.

\subsection{Computational Cost Associated with the Data-Independent Matrix, $\bA_{0}$}
\begin{figure}[h!]
      \centering
      {\parbox{\columnwidth}{
      \centering
      {\includegraphics[width = 0.65\columnwidth, trim=0in 0.9in 0in 0.5in, clip]{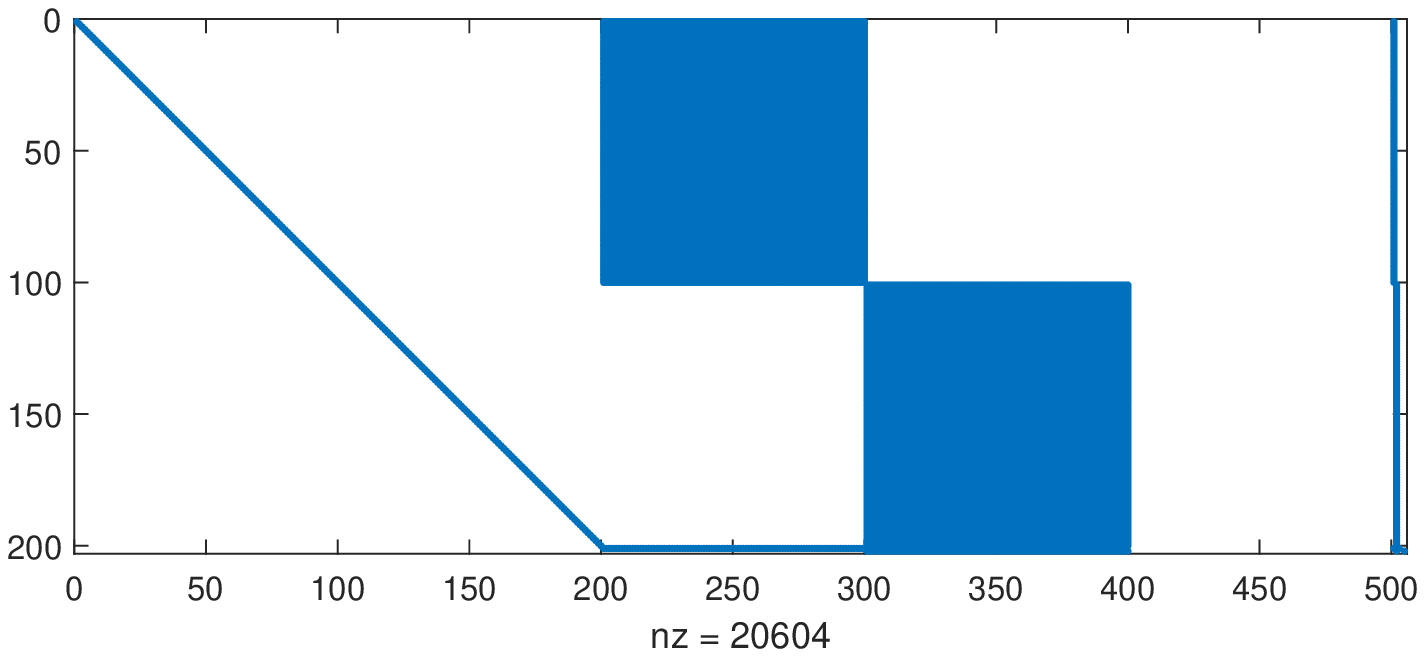}}
      \caption{\textsf{Data-independent sparsity pattern of $\bA_{0}$}}\label{fig:A0}
      }
      }
\end{figure}

The sparsity pattern of $\bA_0$ is shown in Fig.~\ref{fig:A0}.
%
%
The dense blocks in this figure correspond to the $\bB^a$ and $\bB^b$ blocks of $\bA_0$. At first glance, it appears that the Birkhoff matrices might be the main source of a computational burden. To better process this initial intuition, we break up our analysis into two extreme cases of Table~\ref{table:size+speed}, namely rows~1 and 4.
\subsubsection{Thousands of Grid Points}
The first row of Table~\ref{table:size+speed} indicates that storing a Birkhoff matrix for a grid of one thousand points takes only fractions of a gigabyte ($0.008$ GB). Furthermore, solves of \eqref{eq:Ax=b} are milliseconds per iteration even if they are performed at the slowest rate of $O(n^3)$. See top right corner of Table~\ref{table:size+speed}. From the second row of Table~\ref{table:size+speed}, it follows that computations over grids of less than $10,000$ points can also be performed at moderately fast speeds using only a few gigabytes of memory.  In other words, computations over a few thousand grid points are not a serious obstacle to a Birkhoff-theoretic method for trajectory optimization. Coupled with the fact that a universal Birkhoff method has $\mathcal{O}(1)$ condition numbers and a theoretical infinite rate of convergence\cite{newBirk-part-I}, it follows from Table~\ref{table:size+speed} that many aerospace trajectory optimization problems can be solved rapidly and accurately over thousands of grid points.

\subsubsection{Addressing the Million Grid-Point Challenge}
The challenge of solving a million point trajectory optimization problem was posed in \cite{ross-million-2017} as an enabler for advancing a new theory for solving discrete-continuous mission planning problems\cite{ross-TSP-nolcos-2016,ross-TSP-ACC-2019,ross-TSP-arXiv-2020}.  From row~4 of Table~\ref{table:size+speed}, it follows that storing a Birkhoff matrix requires $8$ terrabytes of memory.  Furthermore, solves of \eqref{eq:Ax=b} require a million seconds ($\approx 11$ days) per iteration at $\mathcal{O}(N_n^3)$ speed.  Obviously, these numbers indicate a severe technological barrier.  However, it was shown in \cite{sandia-aas-23} that if the Chebyshev version of the Birkhoff matrix is chosen (see Fig.~\ref{fig:BirkhoffPStypes}) then solves of \eqref{eq:Ax=b} can be performed at $\mathcal{O}(N_n \log N_n)$ computational speed (see column 5 of Table~\ref{table:size+speed}) using matrix-free techniques\cite{matrix-free-ref-1}.  The space complexity of this new approach is $\mathcal{O}(N_n)$ (column 2 of Table~\ref{table:size+speed}).
In other words, the matrix-free, Chebyshev version of the universal Birkhoff-theoretic method holds the potential to solve the million-point trajectory optimization challenge problem\cite{ross-million-2017} in fractions of a second using only a few gigabytes of memory.

\section{A Brief Discussion on Alternative Birkhoff Hessians}\label{sec:NLPHessians}
The Hessian associated with a generic NLP is given by\cite{luenberger-2008,NW:NumOptBook,bazaraa-2006},
\begin{equation}\label{eq:NLP-Hessian}
\partial^2L(\by, \bpsi) := \left[
                               \begin{array}{cc}
                                 L_{y, y} & L_{y, \psi} \\
                                 L_{\psi, y} & 0 \\
                               \end{array}
                             \right]
\end{equation}
where $(\by, \bpsi)$ is the primal-dual pair. Choosing this pair by setting $\by =(\widetilde{X}, \widetilde{V}, \widetilde{U}, x^a, x^b) $,  $\bpsi = (\Lambda, \Omega, \nu, \lambda^b)$ and reusing the special weighted Lagrangian given by \eqref{eq:Lagrangian} we get the following matrix (written in a Hamiltonian form):
\begin{equation}\label{eq:NLPAmat}
\widetilde{\bA} := \left[
             \begin{array}{ccccc;{1pt/2pt}cccc}
              W^{-1}_B H_{X,X} & 0 &W^{-1}_B  H_{X, U}& 0 & 0 & H_{X,\Lambda} & I & 0 & 0  \\
               0 & 0 & 0 & 0 & 0 & -I & \bB^b & 0&  \bb\\ 
               W^{-1}_B H_{U,X}& 0 &  W^{-1}_B H_{U, U}  &  0 & 0& H_{U,\Lambda} & 0 & 0 & 0 \\
               0& 0& 0& \overline{E}_{a,a} & \overline{E}_{a,b}& 0 & -\bw_B^T  &\overline{E}_{a,\nu} & 1 \\
               0&0 & 0 & \overline{E}_{b,a}  & \overline{E}_{b,b}& 0 & 0 & \overline{E}_{b,\nu}  & -1 \\\hdashline[1pt/2pt]
                H_{\Lambda, X} & -I & H_{\Lambda, U} & 0 & 0  & 0 & 0 & 0 & 0  \\
                I &  (\bB^b)^T  & 0 & -\bw_B& 0 & 0 & 0 & 0 & 0 \\
                0 & 0 & 0 & \overline{E}_{\nu,a} & \overline{E}_{\nu,b}& 0 &0 & 0 & 0\\
                0 & \bb^T & 0 & 1 & -1 &0 &0 &0 & 0 \\
             \end{array}
           \right]
\end{equation}
Note the following:
\begin{enumerate}
\item Unlike $\bA$ (Cf.~\eqref{eq:Hessian=}) which is a ``$10 \times 10$'' matrix, $\widetilde{\bA}$ is  ``$9 \times 9$''.  This is because $\lambda^a$ is absent in \eqref{eq:NLPAmat}. As indicated in Remark~\ref{rem:lam-a}, \eqref{eq:lam-a-def}, and hence $\lambda^a$, does not exist in the context of an NLP. Instead, \eqref{eq:tvc-a} takes the following form,
$$  \lambda^b -  \bw_B^T \Omega + \partial_{x^a} \overline{E}(\nu, x^a, x^b) =  0$$
\item  Unlike $\bA$ which is not symmetric, $\widetilde{\bA}$ is indeed symmetric.
\item The bottom left block of $\widetilde{\bA}$ is not the Jacobian of the constraints as is the case with an unweighted Lagrangian.
\item A direct application of Gershgorin's theorem to $\widetilde{\bA}$ will not yield Theorem~\ref{thm:ross-hess}.
\item An application of Gershgorin's theorem to $\widetilde{\bA}$ will yield eigenvalue estimates that are grid dependent.  See, for example, the last column/row.
\end{enumerate}

Now suppose a standard unweighted Lagrangian is used instead of the special one given by \eqref{eq:Lagrangian}. Retracing all the steps taken in Sections~\ref{sec:Lag2Ham} and \ref{sec:HamilFormOfHess}, it can be observed (see also \cite{newBirk-part-I} for additional details) that \eqref{eq:gradL=0-pre} generates a weak implementation of the Hamiltonian conditions indicated in \eqref{eq:gradL2H=0}.  To better appreciate why a weak implementation may not be desirable, consider its implications on \eqref{eq:Hu=0-strong}.  Its weak implementation would take the form,
\begin{equation}\label{eq:Hu=0-weak}
 \bw_B \circ \partial_u H(\Lambda, X, U)  =  0
\end{equation}
Because it can only be imposed up to a practical tolerance of $\delta > 0$, the computational version of \eqref{eq:Hu=0-weak} takes the form,
\begin{equation}\label{eq:Hu=0-weak-delta}
 \bw_B \circ \partial_u H(\Lambda, X, U)  \le  \delta [1, \ldots, 1]^T
\end{equation}
Hence the dual computational constraint that effectively gets imposed by \eqref{eq:Hu=0-weak-delta} can be written as,
\begin{equation}\label{eq:Hu=d/w}
\partial_u H(\lambda_i, x_i, u_i)  \sim  \mathcal{O}(\delta/w_i), i = 0, \ldots, N
\end{equation}
Because $ 0< w_i < 1, i = 0, \ldots, N$, the effective tolerance on imposing Pontryagin's Hamiltonian minimization condition is greater than $\delta$, with the ``worst'' impositions occurring at the extreme points, i.e., $0$ and $N$.  This is because $w_0$ and $w_N$ are the smallest of all the numbers $ 0< w_i < 1, i = 0, \ldots, N$\cite{newBirk-2023}.  In other words, $u_i, i = 0, \ldots, N$ resulting from solving \eqref{eq:Hu=d/w} would be inaccurate (relative to an imposition of \eqref{eq:Hu=0-strong}) with the worst accuracy near the boundary points. For these reasons, a standard Lagrangian is not computationally efficient, particularly for increasing values of $N$.  The alternative is, of course, to use a weighted Lagrangian.

We briefly note that the argument for a weighted Lagrangian is not limited to Birkhoff-theoretic methods; it has been discussed before\cite{advances,ross:shooting-arXiv} in the context of Lagrange pseudospectral methods (i.e., pseudospectral methods based on differentiation matrices\cite{advances}) and Runge-Kutta methods\cite{ross:shooting-arXiv}.

As a final point of discussion, consider once again retracing all the steps taken in Sections~\ref{sec:Lag2Ham} and \ref{sec:HamilFormOfHess}. It can be easily observed by this process that the source of asymmetry in $\bA$ is the introduction of \eqref{eq:lam-a-def} and the use of the identity,
\begin{equation}\label{eq:b=identity}
\bb \equiv W_B^{-1} \bw_B
\end{equation}
Both \eqref{eq:lam-a-def} and \eqref{eq:b=identity} are necessary to cast the Lagrange form of the necessary conditions (i.e., \eqref{eq:gradL=0-pre}) to an ``equivalent'' Pontryagin Hamiltonian form given by \eqref{eq:gradL2H=0}.
In choosing to write a Hessian  in its strong Hamiltonian form with both initial and final transversality conditions, the symmetry of \eqref{eq:NLPAmat} is lost in \eqref{eq:Hessian=}.  Nevertheless, note that the asymmetric and strong Hessian given by \eqref{eq:Hessian=} is more meaningful in directly connecting it to Pontryagin's Principle\cite{ross-book} in at least two ways: (i) the explicit use of both the initial and final transversality conditions and (ii) a direct connection to the Pontryagin Hamiltonian without an intervening weight function; see also \eqref{eq:PA=}.

Because the entire body of literature in NLPs is associated with a symmetric Hessian, it is apparent that new algorithms beyond NLPs are needed to take full advantage of the properties associated with \eqref{eq:Hessian=} such as those identified in Theorem~\ref{thm:ross-hess}.  This is an open area of research.

\section{Conclusions}
In much of the literature on trajectory optimization, it is more common to see a detailed analysis of the Jacobian and its sparsity pattern rather than the Hessian. This is likely because in typical (non-Birkhoff) discretizations of optimal control problems, a detailed analysis of the Hessian is quite foreboding.  For example, in Runge-Kutta methods, functions must be evaluated at collocation points that are not the discretization points.  In sharp contrast, a universal Birkhoff-theoretic method decouples the discretization from the problem data functions through a natural introduction of virtual and co-virtual control variables.  These primal and dual virtual variables appear linearly in the discretized optimal control problem and hence facilitate a substantially simpler process of analyzing a Hessian at a detailed granular level through the use of the covector mapping principle. The main ``trick'' employed to standardize the analysis in this note is to utilize the discretized Pontryagin Hamiltonian and its derivatives. Furthermore, the nonlinear components of the Hessian do not contain any Birkhoff-specific elements.  This makes a Birkhoff Hessian more amenable to a detailed analysis. All of these desirable features point toward additional advantageous attributes of the universal Birkhoff theory that go well beyond its previously-discussed properties of flat condition numbers. Obviously, much remains to be done in advancing the theory for fast trajectory optimization.  Nonetheless, the results obtained so far lend further credibility to the notion that the universal Birkhoff theory offers a disruptive method for rapidly solving nonlinear, nonconvex trajectory optimization problems.

\section*{Acknowledgments}
The author gratefully acknowledges partial funding for this research provided by the Air Force Office of Scientific Research (AFOSR) and the Defense Advanced Research Project Agency (DARPA).  The views and conclusions contained herein are those of the author and should not be interpreted as necessarily representing the official policies or endorsements, either expressed or implied, of the AFOSR, or DARPA or the U.S. Government.




\end{document}